# Modeling the SBC Tanzania Production-Distribution Logistics Network

Nelson Christopher Dzupire[1*], Yaw Nkansah-Gyekye [1], Silas Steven Mirau[1]

[1]Nelson Mandela African Institution of Science and Engineering,

School of Computational and Communication Science and Engineering,

P. O. Box 447, Arusha, Tanzania.

* E-mail of corresponding author: dzupiren@nm-aist.ac.tz

**Abstract**

The increase in customer expectation in terms of cost and services rendered, coupled with competitive business environment and uncertainty in cost of raw materials have posed challenges on effective supply chain engineering  making it essential to do cost-benefit analysis before making final decisions on production-distribution logistics. This paper provides a conceptual model that provide guidance in supply chain decision making for business expansion. It presents a mathematical model for production-distribution of an integrated supply chain derived from current operations of SBC Tanzania Ltd which is a major supply chain that manages products' distribution in whole of Tanzania. In addition to finding the optimal cost, we also carried out a sensitivity analysis on the model so as to find ways in which the company can expand at optimal cost, while meeting customers' demands. Genetic algorithms is used to run the simulation for their efficient in solving combinatorial problems.

**1. Introduction**

Nowadays the success of a company depends mostly on management of its supply chain. The supply chain management (SCM) is defined as a set of approaches used to effectively integrate suppliers, manufacturers, distribution centers so that goods are produced and distributed at the right quantities to the right locations, and at the right time in order to minimize system wide cost while satisfying customer service level requirements (Kuo and Han, 2011). It involves the integration of business processes from customers through the suppliers that provide products, services and information.

In general the complexity of the business environment has been challenged by several factors namely expansion of the market, wide range of suppliers, increased competition, and customer's demands on the performance of a company in terms of waiting time, cost and quality of the product (Copacino, 1997). This has brought several questions like where to best site warehouses (distribution centers) and manufacturing plants, number of plants and warehouses to work with in an industry and their corresponding capacities. Cheng and Lin (2008) stated that typical SCM goals include transportation network design, plant and distribution center location and allocation, production schedule streamlining and efforts to improve order response time. This has made logistics network design a comprehensive strategic decision problem that has to be optimized for long term benefit of the whole supply chain.

Production-distribution logistics model designing has attracted the attention of many researchers (Mehdizadeh and Afrabandpei, 2012; Amiri, 2006; Xiaobo *et al.*, 2007). Through it, we determine the number, location, capacity and type of plants  and distribution centers to be used. We can also design the distribution system to several centers and retailers, and the amount of raw materials to consume, quantities of products to produce. So production-distribution logistics of a supply chain covers the entire process of buying and transporting raw materials to plants, conversion of these raw materials into products, the transportation of the products  to various distribution centers and eventual delivery to retailers. As such production-distribution logistics management reduces total costs and hence it is a key issue in today's competitive business environment especially multinational companies (Cheng and Lin, 2008).

In this research we develop a production-distribution logistics model of SBC Tanzania Ltd which is a supply chain that produces soft drinks. It was established in Tanzania in 2001 and produces five products namely milinda, mountain dew, pepsi, 7up and evervess in different flavors and volumes (Sbc Tanzania, 2010). They have 200 wholesaler and 5220 retailer customers currently. The model is designed by considering minimization of cost involved in the production and distribution process with sensitivity analysis so that we investigate the model performance and illustrate how parameter changes influences the feasibility of the business. We will consider the following parameter analysis: varying distance coverage to reach retailers, varying distribution center's storage capacity, and varying number of plants. As observed by Farahani and Elahipanah (2008) key to the success of any business is satisfying customer's demands on time which may result in cost reductions and increase in service level, hence we would like to find out if variation of any parameter in the model satisfies the customers' demands and is business viable for an industry. The analysis is based on the integration of several functions of the supply chain into a single optimization model, in such a way that we simultaneously optimize several decision variables of different cost functions.

Shen (2007) in a recent survey observed that most of the problems in supply chain optimization are combinatorial and NP-hard, and therefore difficult to solve. The logistics network design problem can as well be categorized as combinatorial and NP-hard, of which several approaches of solving such problems exist, classified as exact, heuristics and meta-heuristic (Mehdizadeh and Afrabandpei, 2012). Recently genetic algorithms have received overwhelming attention as an approach to solve such optimization problems. The genetic algorithms are known to be efficient-solving and easy adaptive, especially where traditional methods failed to provide good solutions (Liao and Hsieh, 2009; Pinto, 2004). They have been used in optimizing logistics network problems and different ways of chromosome representation are proposed depending on specific problem.

A production-distribution logistic network design problem involves simultaneous decisions about location of plants and distribution centers, specification of plants and distribution centers' capacities and distribution systems for raw materials and products in-process. A great deal of research has been carried out to develop mathematical models that represent such decisions and the integration of the supply chain design, with exact or heuristic algorithms have been used to solve such models. Park* (2005) proposed a solution for an integrated production and distribution planning and investigated the effectiveness of the integration in a multi-plant, multi-retailer logistic environment, with maximizing profit as the objective function. Using heuristic algorithm, results were computed which confirmed substantial advantage of the integrated planning approach over decoupled ones. Furthermore, parameter analysis indicated that with right conditions, it is effective to integrate production and distribution functions. In the work of Lejeune (2006) a mixed integer programming model with cost minimization as an objective function in a three echelon supply chain network was developed. The model had a sustainable inventory production distribution plan which was constructed over a multi-period horizon. Results were computed using variable neighborhood algorithm. Mehdizadeh and Afrabandpei (2012) designed a multi-stage and multi-product logistic network model which was a mixed nonlinear integer programming model. The objective function was minimizing transport and holding cost in a three echelon supply chain.

This paper is organized as follows: section 2 gives the problem description in which we introduce the objective function and its assumptions followed by the mathematical models for the integrated supply chain. The proposed methodology is described in section 3 followed by computed results and discussion in section 4. Finally, the conclusion is given in section 5.

## 2.1 Problem Description

The major purpose of this work is to investigate effectiveness and feasibility of an integrated production-distribution logistic model of SBC Tanzania Ltd by varying the following parameters: number of plants, capacities of distribution centers and reducing distance covered by retailers. The supply chain involves five suppliers some of which are from outside Tanzania, four manufacturing plants, four distribution centers and many customers namely 5220 retailers and 200 wholesalers.

The objective function is to minimize total cost which involves the following: transportation cost of raw materials from suppliers to plants and the cost of buying the raw materials, transportation cost of delivering the products to distribution centers and to retailers, and holding cost of products at distribution centers. As a company the aim is to minimize the cost as much as possible while maintaining quality and efficiency, and at the same time meeting customers' demands on time so that profitability of the entire supply chain is maximized. In addition they also aim to reach out to as many customers as possible which may result in the expansion of the company. Here are the assumptions for the model to be formulated:

1. Number of retailers and suppliers and the capacity of the suppliers are known.
2. Number of plants and distribution centers, and their capacities are also known.
3. Demands of customers are uncertain but can be determined from past history.
4. Each plant receives raw materials from all the suppliers.
5. Each retailer is served from a single distribution center whereas each distribution centers gets products from all the plants.

## 2.2 Problem Formulation

A mixed integer mathematical programming model for the production-distribution logistics design of an integrated supply chain is presented. In presenting the model we use the following:

    a) Notation for indices of the entities are as follows:

*s*: Suppliers

*k*: Plants

*j*: Distribution centers

*i*: Retailers

*t*: time

b) Variables for quantities are as follows:

$r_{sk}$: quantity of raw materials from supplier *s* to plant *k*,

$p_{kj}$: quantity of products manufactured at plant *k* delivered to distribution center *j*,

$t_{ji}$: quantity of products transported from distribution center *j* to retailer *i*.

c) The variable notation for model parameters are:

$D_k$: capacity of plant *k*,

$C_s$: capacity of supplier *s*,

$H_j$: holding capacity of distribution center *j*,

$DC_{tot}$: total number of distribution centers in the supply chain,

$P_{tot}$: total number of manufacturing plants in the supply chain,

$d_i$: quantity of demanded products at retailer *i*,

$c_s$: unit cost of buying and transporting raw materials at supplier *s*,

$h_j$: unit storage cost of products at distribution center *j*,

$c_{kj}$: unit transportation cost of products from plant *k* to distribution center *j*,

$r_{ji}$: unit transportation cost of products from distribution center *j* to retailer *i*.

The mathematical formulation of the model is as follows: minimize

$$\sum_s \sum_k c_s r_{sk} + \sum_k \sum_j c_{kj} p_{kj} + \sum_t \sum_j h_{jt} p_{kj} + \sum_j \sum_i r_{ji} t_{ji} \qquad (1.1)$$

Subject to the following constraints:

$$\sum_i d_i \leq \sum_j H_j \qquad (1.2)$$

$$\sum_k p_{kj} \geq \sum_i q_{ji} \qquad (1.3)$$

$$\sum_j t_{ji} = d_i \qquad (1.4)$$

$$u \sum_j p_{kj} \leq \sum_s r_{sk} \qquad (1.5)$$

$$u \sum_j p_{kj} \leq D_k \qquad (1.6)$$

$$\sum_k r_{sk} \leq C_s \qquad (1.7)$$

$$p_{kj}, r_{sk}, c_{kj}, H_j, C_s, u, t_{ij}, h_j, D_k, d_i, r_{ji}, c_s \geq 0 \qquad (1.8)$$

The objective function (1.1) expresses the total cost of the supply chain. The total cost involves buying raw materials, storage and transportation of raw materials to plants and products to distribution centers and then to retailers. Constraint (1.2) implies that all demanded products should not exceed the storage capacity of all distribution centers. Constraint (1.3) limits amount of products sent to retailers to be with the range of those that are produced by manufacturing plants. Customer satisfaction is guaranteed in constraint (1.4) where all demands are met without backlog. Constraint (1.5) limits quantity of products produced to the available quantity of raw materials. Similarly quantity of products produced should within the capacity of plants (1.6) where we have *u* as the utilization factor. Finally, constraint (1.7) implies that amount of raw materials sent to manufacturing plants is within the capacity of the suppliers with constraint (1.8) ensuring positivity.

## 3. Genetic Algorithms

With optimization problems that are combinatorial and NP-hard, it is hard to use algorithms that can find exact optimal solutions. Even though some such algorithms exist like branch-and-bound, branch-and-cut and branch-and-prize, but with real size problems it is time consuming and have higher computational complexity (Mehdizadeh and Afrabandpei, 2012; Zitzler *et al.*, 2004; Konak *et al.*, 2006). The exact algorithms solve the mixed integer program problems by progressing from node to node to implicitly exhaust all possible combinations. This makes the algorithm infeasible as the number of combinations grow exponentially with the size of the problem (Lejeune, 2006). Therefore we use algorithms that reasonably approximate the solution in polynomial time. One such algorithm is the genetic algorithms which is one of the evolutionary algorithms. The concept of genetic algorithms was developed by Holland (1975) and his colleagues. It was inspired by the evolution theory whereby weak and unfit species in their environment goes to extinction by natural selection while the strong ones survive and pass on their genes to the future generation through reproduction. A solution in genetic algorithm is referred as an individual or chromosome and it is made up of genes which controls one or more features in the chromosome. These genes are either binary or real coded when implementing depending on the size of the problem and user preferences. A collection of chromosomes forms a population. Normally initial population is randomly generated but as time goes by it includes fitter and fitter chromosomes until it converges to a single solution. The algorithm works well on mixed combinatorial problems whether it is discrete or continuous and is less susceptible to converge to a local optimal as compared to gradient search methods (Danalakshmi and Kumar, 2012). The value of the chromosome as calculated using the objective function is called its fitness value. Genetic algorithms uses two operators namely crossover and mutation to generate a new population from existing the existing population. With crossover, two chromosomes, referred as parents combine together to form new chromosomes called offspring. The parents selected are those with better fitness values. The mutation operator introduces random changes into characteristics of a chromosome by changing part of it. Mutation operator enables diversity of the population so as not to converge to a local optima solution. Then we have a reproduction process that involves selection of chromosomes to form the next generation. In general the fitness of each chromosome determines its probability of survival to the next generation. Several selection approaches are available depending on how the fitness value is calculated and include the following: proportional, ranking, and tournament selection. In each generation, fitness value of every chromosome is evaluated, multiple chromosomes are stochastically selected from the current population based on their fitness, and are recombined to form a new population. The algorithm terminates when either the maximum number of generations is reached or a satisfactory fitness level has been reached.

In this research we used the fast non-dominated sorting genetic algorithm (NSGA-II) developed by Deb *et al.* (2002) which has proved to be quiet efficient in many applications and its performance is far much better than most existing ones. NSGA-II is an elitist and fast strategy, modular and flexible, emphases on the non dominated solutions, can be applied to a large wide of problems and can easily be implemented in the global optimization toolbox in MATLAB. It is an improvement of the non-dominated sorting genetic algorithm.

Without loss of generality, considering a minimization problem, a solution A is said to dominate solution B if A has a lower value for at least one of the objective functions and is not worse than B in the remaining objective functions. So a solution is non-dominated if no solution dominates it. In this algorithm, for each individual solution we calculate the number of solutions that dominate it, and determine a set of solutions that this solution dominate. Thus all solutions are ranked into non-dominated fronts based on their ranks calculated by how many solutions dominate them. Hence those on rank 1 are best and non-dominated. To allow diversification, we also compute the crowding distance of each solution so that we find the solution density surrounding a particular solution in the population. To enable elitism, individuals for the next generation are selected from both the parents and children based on their non-domination front and crowding distance, starting from those with lower ranks, and for those in the same rank, we prefer those with higher crowding distance, until N individuals are selected (Deb *et al.*, 2002; Konak *et al.*, 2006; Zitzler *et al.*, 2004).

## 4. Results and Discussion

We simulated three scenarios in this paper as follows: the current state, changing the capacities of distribution centers, changing number of both plants and distribution centers and their capacities. The simulation was done in MATLAB 2013a with the following parameters: crossover probability of 0.6, mutation of 0.001 and a population size of 50. We also used ranking as a scaling function. The initial population was uniformly generated.

4.1 Current State:

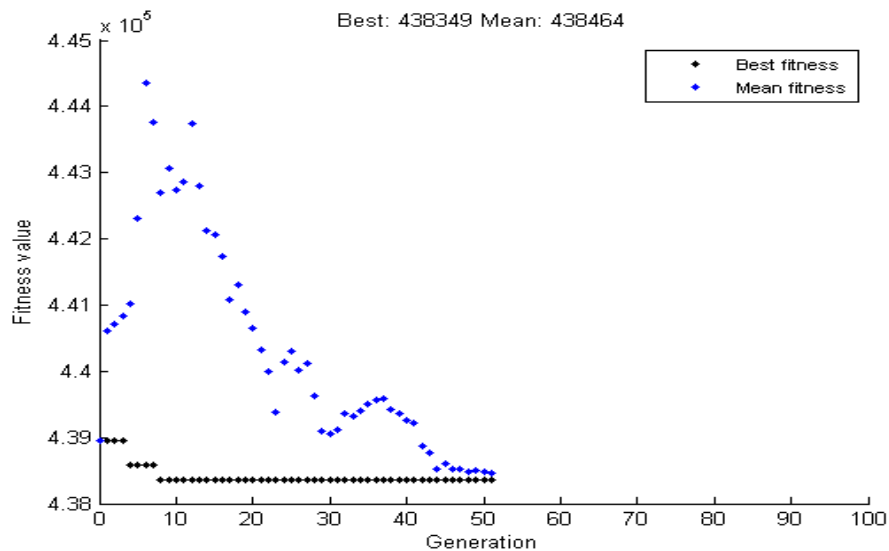

Fig 1. Minimum Cost ('00' TZ Shillings) at the current state

With the current state of the SBC Tanzania supply chain, the optimal cost of production-distribution logistics model is found to be 43,834,900 TZ shillings per week. The amount of products produced per day in cases, with each case containing 24 bottles, is as follows:

Table 1. Number of cases

|  | Plant 1 | Plant 2 | Plant 3 | Plant 4 |
| --- | --- | --- | --- | --- |
| Distribution center 1 | 315 | 0.0 | 0.0 | 11196 |
| Distribution center 2 | 0.0 | 4949 | 8964 | 0.0 |
| Distribution center 3 | 10474 | 0.0 | 0.0 | 1897 |
| Distribution center 4 | 0.0 | 5932 | 6766 | 0.0 |

There are currently four plants with capacities 12800, 12000, 25600, and 12800 cases respectively. There are also four distribution centers with a capacity of 12000 cases each. In this model, we found out that it is optimal if two plants can produce all the products between them so that each distribution center can be supplied by only two manufacturing plants. This will also be in tandem with reducing distance coverage to distribution to centers.

4.2 Changing Storage Capacities of Distribution Centers:

In this analysis, we changed the storage capacity of distribution centers from the current 12000 cases to 15,000 cases for each. In this way we tried to create a lot of storage space for all the four centers so that we are able to meet customer's demands on time. Below is the graph for the analysis and its corresponding optimal cost.

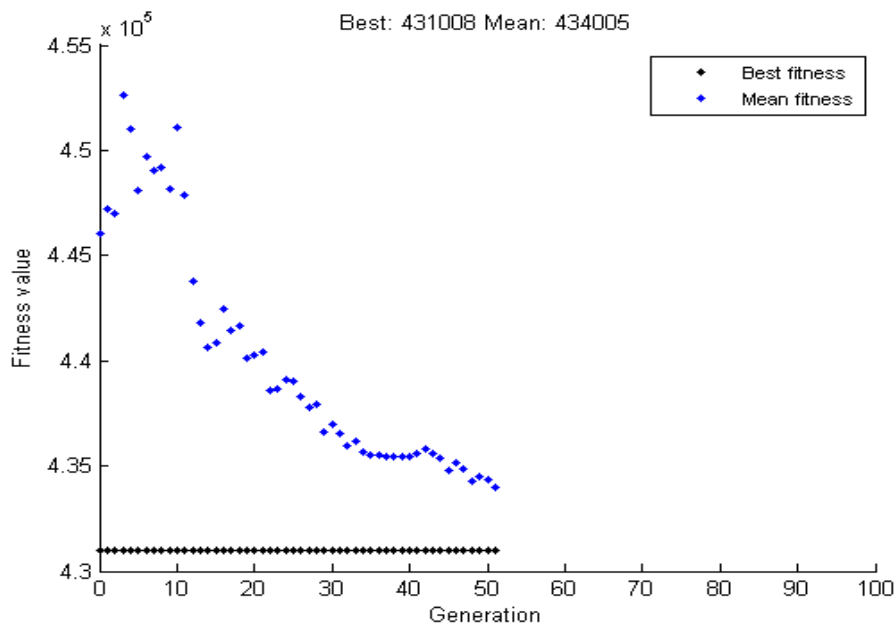

Figure 2. Optimal cost in ('00'TZ shillings)

The production schedule per day (in cases) with this arrangement is as follows:

Table 2. Number of cases with changed distribution centers' capacity

|  | Plant 1 | Plant 2 | Plant 3 | Plant 4 |
|---|---|---|---|---|
| Distribution center 1 | 4257 | 0.0 | 0 | 10253 |
| Distribution center 2 | 7363 | 0.0 | 745 | 6363 |
| Distribution center 3 | 1670 | 8753 | 0.0 | 4572 |
| Distribution center 4 | 0.0 | 3427 | 11155 | 0.0 |

In this case, the company would be spending 43,100,800 TZ shillings which is slightly lower as compared to the previous current state. We observe that some distribution centers are receiving products from at least three plants which is not the case in the other case. The reduction in cost can be explained in that we still worked with the same number of customers and their demands in the previous case, however if the demand for products increases, then the cost of the production-distribution model may increase slightly.

4.3 Varying Distance Coverage to Reach Retailers and Inventory cost:

In this case we simulate increasing the number of both plants and distribution centers so that they are widely distributed in the country hence easy to reach retailers. It is assumed that increasing the number of plants and distribution center will enable the company to reach to as many customers as possible and at the same time reduces the distance travelled to reach retailers so that we can be able to vary inventory cost as transportation cost is reduced. In addition, we also increase the capacities of both plants and distribution centers. Of the four plants, three plants have capacities less than 15,000 cases so we increase it to 15,000 cases and the other one currently at 25,600 cases to 30,000 cases hence we have six plants with capacity of 15,000 and one plant with 30,000 capacity in total as seven plants. All distribution centers have a capacity of 12,000 cases so we set it to 15,000 cases and increase their number to eight. However the simulation does not involve demands by customers so the only limitation is the capacities of both plants and distribution centers. We got the following graph and its corresponding daily production-distribution table below:

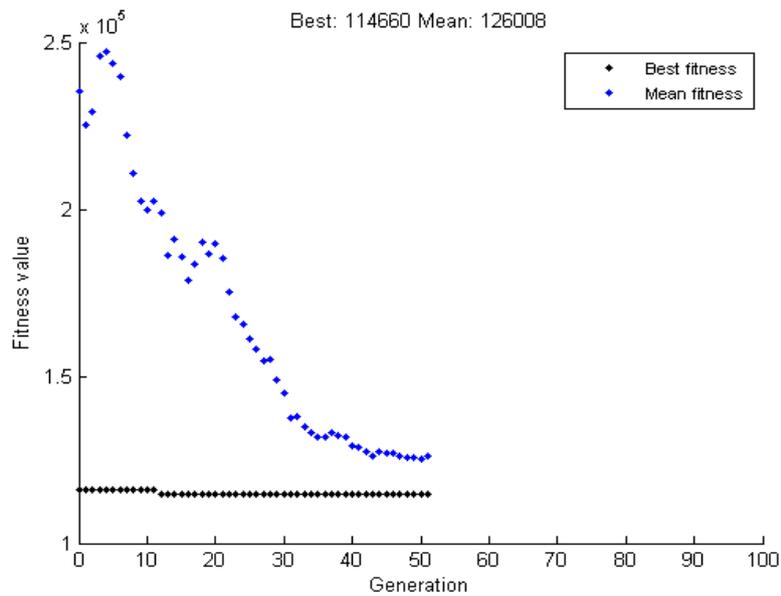

Fig 3. Optimal Cost ('000') TZ shillings

Table 3. Daily schedule of production (in cases)

|  | Plant 1 | Plant 2 | Plant 3 | Plant 4 | Plant 5 | Plant 6 | Plant 7 |
|---|---|---|---|---|---|---|---|
| DC 1 | 800 | 0.0 | 0.0 | 11364 | 0.0 | 2452 | 0.0 |
| DC 2 | 0.0 | 0.0 | 4215 | 0.0 | 6021 | 0.0 | 4750 |
| DC 3 | 5688 | 3061 | 4535 | 1074 | 0.0 | 0.0 | 0.0 |
| DC 5 | 2133 | 6438 | 800 | 4957 | 0.0 | 0.0 | 0.0 |
| DC 4 | 0.0 | 0.0 | 0.0 | 5800 | 0.0 | 2944 | 6167 |
| DC 6 | 0.0 | 4990 | 644 | 0.0 | 5717 | 0.0 | 3482 |
| DC 7 | 0.0 | 0.0 | 0.0 | 6632 | 0.0 | 7452 | 0.0 |
| DC 8 | 6128 | 0.0 | 4086 | 0.0 | 3111 | 1669 | 0.0 |

DC: Distribution center

In this case we have an optimal cost of 114,660,000 TZ Shillings per week to produce a daily production as presented in Table 3. Worth mentioning is that production distribution was done with no demand constraints and the only limitations were capacities of both plants and distribution centers. In this we have a lot of products being produced and the company may benefit if there can be an increase in demands and widely distribute the manufacturing plants and distribution centers so that they reach out to several customers.

## 5. Conclusion

In all the cases we have considered here, the major difference is the third case where the optimal amount is 114,660,000 TZ Shillings. In the first case where the optimal cost is 43,834,900 TZ Shillings, we discover that it is small as compared to the cost being incurred now which ranges from 45-65 million TZ shillings per week meeting the same demands. So we urge the company to explore our production-distribution logistics as it is feasible and cheaper. In the second case, the lower cost being found compared to the first case can be compensated with fixed cost of maintaining those big distribution centers whose capacities have been increased. Though it is feasible and cheaper, the company may have to decide whether buying or establishing new distribution centers with improved capacities is viable since it only offers a 13% increase in number of products produced compared to the first one. So this model offers the company an opportunity to expand and reach out to many customers. The last case offers much expansion for the company, as its production schedule represents 57% increase from the current model. The production-distribution schedule is far better as compared to the two cases discussed. Its feasibility depends on two fronts: company having money to establish the manufacturing plants and distribution centers in addition to meeting the cost of production per week, and also if there can be demands so that products produced does not end up overcrowded in the distribution centers. However it allows the company to reach out to many customers and have distance travelled to reach retailers reduced, such that they can be able to vary inventory cost. So if the company can establish market for all the products then they can explore this case. It should be noted that the choice of the variation in both capacities and number of plants or distribution centers was done arbitrary after exploring several cases.

Future research should explore a scenario where the company has a specific amount of money to expand the business. So in such case, researchers should advise the company what would be an optimal production-distribution schedule using such an amount of money. In that way the company will be able to expand gradually other than at one go.

**Acknowledgement**

Nelson is indebted to the Nelson Mandela African Institution of Science and Technology for financial and material support.

**References**

Amiri, A. (2006). Designing a distribution network in a supply chain system: Formulation and efficient solution procedure. *European Journal of Operational Research*. **171**: 567-576.

Cheng, R. and Lin, L. (2008). Network models and optimization: Multiobjective genetic algorithm approach. Springerpp.

Copacino, W.C. (1997). Supply chain management: The basics and beyond. CRC Presspp.

Danalakshmi, C.S. and Kumar, G.M. (2012). Optimization of Supply Chain Network Using Genetic Algorithm. *J. Manuf. g Eng*.

Deb, K., Pratap, A., Agarwal, S. and Meyarivan, T. (2002). A fast and elitist multiobjective genetic algorithm: NSGA-II. *Evolutionary Computation, IEEE Transactions on*. **6**: 182-197.

Farahani, R.Z. and Elahipanah, M. (2008). A genetic algorithm to optimize the total cost and service level for just-in-time distribution in a supply chain. *International Journal of Production Economics*. **111**: 229-243.

Holland, J.H. (1975). Adaptation in natural and artificial systems: An introductory analysis with applications to biology, control, and artificial intelligence. U Michigan Presspp.

Konak, A., Coit, D.W. and Smith, A.E. (2006). Multi-objective optimization using genetic algorithms: A tutorial. *Reliability Engineering & System Safety*. **91**: 992-1007.

Kuo, R. and Han, Y. (2011). A hybrid of genetic algorithm and particle swarm optimization for solving bi-level linear programming problem–A case study on supply chain model. *Applied Mathematical Modelling*. **35**: 3905-3917.

Lejeune, M.A. (2006). A variable neighborhood decomposition search method for supply chain management planning problems. *European Journal of Operational Research*. **175**: 959-976.

Liao, S.-H. and Hsieh, C.-L. (Ed.) (2009). An integrated inventory control and facility location system with capacity constraints: A multi-objective evolutionary approach. *In: 20th Annual Conf. Production and Operations Management Society, Orlando, Florida*, pp. 1-4, 2009.

Mehdizadeh, E. and Afrabandpei, F. (2012). Design of a Mathematical Model for Logistic Network in a Multi-Stage Multi-Product Supply Chain Network and Developing a Metaheuristic Algorithm. *Journal of Optimization in Industrial Engineering*. **10**: 35-43.

Park*, Y. (2005). An integrated approach for production and distribution planning in supply chain management. *International Journal of Production Research*. **43**: 1205-1224.

Pinto, E.G. (2004). Supply chain optimization using multi-objective evolutionary algorithms. *Retrieved December*. **15**: 2004.

SBC Tanzania (2010). http://www.sbctanzania.com/company-profile.html.

Shen, Z. (2007). Integrated supply chain design models: a survey and future research directions. *Journal of Industrial and Management Optimization*. **3**: 1.

Xiaobo, Z., Xu, D., Zhang, H. and He, Q.-M. (2007). Modeling and analysis of a supply–assembly–store chain. *European Journal of Operational Research*. **176**: 275-294.

Zitzler, E., Laumanns, M. and Bleuler, S. (2004). A tutorial on evolutionary multiobjective optimization, Metaheuristics for multiobjective optimisation. Springer, pp. 3-37.